\documentclass[12pt]{amsart}
\usepackage{amssymb}
\usepackage[all]{xy}
\usepackage[pdfencoding=auto]{hyperref}
\usepackage{fullpage}
\allowdisplaybreaks

\numberwithin{equation}{section}

\theoremstyle{plain}
\newtheorem{theorem}[equation]{Theorem}
\newtheorem{lemma}[equation]{Lemma}
\newtheorem{prop}[equation]{Proposition}
\newtheorem{cor}[equation]{Corollary}

\theoremstyle{definition}
\newtheorem{defn}[equation]{Definition}

\newtheorem{chunk}[equation]{}

\theoremstyle{remark} 
\newtheorem{remark}[equation]{Remark} 
\newtheorem{rk}[equation]{Remark}

\newcommand{\bF}{\mathbb{F}}
\newcommand{\ii}{\mathsf{i}\,}

\newcommand{\bZ}{\mathbb{Z}}
\newcommand{\bfA}{\mathbf{A}}
\newcommand{\cB}{\mathcal{B}}
\newcommand{\cG}{\mathcal{G}}
\newcommand{\cO}{\mathcal{O}}
\newcommand{\cP}{\mathcal{P}}
\newcommand{\cU}{\mathcal{U}}
\newcommand{\fA}{\mathfrak{A}}
\newcommand{\fM}{\mathfrak{M}}

\newcommand{\Aut}{\mathsf{Aut}}

\renewcommand{\ge}{\geqslant}
\renewcommand{\geq}{\geqslant}
\newcommand{\gr}{\mathsf{gr}}
\newcommand{\Hom}{\mathsf{Hom}}
\newcommand{\Irr}{\mathsf{Irr}}
\newcommand{\Jen}{\mathsf{Jen}}
\newcommand{\Ker}{\mathsf{Ker}}
\renewcommand{\le}{\leqslant}
\renewcommand{\leq}{\leqslant}
\newcommand{\Mat}{\mathsf{Mat}}

\newcommand{\w}{{\sf w}}
\newcommand{\x}{{\sf x}}
\newcommand{\y}{{\sf y}}
\newcommand{\z}{{\sf z}}

\author{David Benson, Radha Kessar, and Markus Linckelmann}

\address{David Benson \\ 
Institute of Mathematics\\ 
Fraser Noble Building\\
University of Aberdeen\\ 
King's College\\ 
Aberdeen AB24 3UE\\ 
United Kingdom}

\address{Radha Kessar and Markus Linckelmann \\
School of Mathematics, Computer Science \& Engineering \\
Department of Mathematics \\
City, University of London \\
Northampton Square \\
London EC1V 0HB \\
United Kingdom}



\title{Structure of blocks with normal defect and abelian $p'$ inertial quotient}

\begin{document}

\begin{abstract}
Let $k$ be an algebraically closed field of prime characteristic $p$.
Let $kGe$ be a block of a group algebra of a finite group  $G$, with
normal defect group $P$ and abelian $p'$ inertial quotient $L$.
Then we show that $kGe$ is a matrix algebra over a quantised
version of the group algebra of a semidirect product of $P$ with
a certain subgroup of $L$. To do this, we first examine the
associated graded algebra, using a Jennings--Quillen style
theorem.

As an example, we calculate the
associated graded of the basic algebra of the non-principal
block in the case of a semidirect product
of an extraspecial $p$-group $P$ of exponent $p$ and order $p^3$ with
a quaternion group of order eight with the centre acting trivially.
In the case $p=3$ we give explicit generators and relations for the
basic algebra as a quantised version of $kP$.
As a second example, we give explicit
generators and relations in the case of a group 
of shape $2^{1+4}:3^{1+2}$ in characteristic two.
\end{abstract}

\maketitle

\section{Introduction} 

Throughout this paper $p$ is a prime and $k$ is an algebraically closed field of 
characteristic $p$.
The study of blocks with normal defect groups has a long history,
starting with the work of Brauer~\cite{Brauer:1956a}, and continuing
with Reynolds~\cite{Reynolds:1963a}, Dade~\cite{Dade:1973a}, and
K\"ulshammer~\cite{Kulshammer:1985a}. In the case of abelian normal
defect, abelian inertial quotient and one simple module, explicit
descriptions of the basic algebra were given by 
Benson and Green~\cite{Benson/Green:2004a}, and 
Holloway and Kessar~\cite{Holloway/Kessar:2005a}.  Dropping 
the hypothesis of one simple module led to our 
paper~\cite{Benson/Kessar/Linckelmann:2019a}.  The main structural feature of the
basic algebras calculated in these papers is that they appear to be quantised
versions of the group algebras of semidirect products of a defect group and
a subgroup of the inertial quotient.

The purpose of this paper is to generalise the results from 
\cite{Benson/Kessar/Linckelmann:2019a} to 
 blocks of group algebras over $k$  of finite groups that have a normal defect group $P$ 
 which is no longer necessarily abelian,  but still with abelian $p'$ inertial
quotient $L$.
By a theorem of K\"ulshammer~\cite{Kulshammer:1985a}, 
any such block is isomorphic to a matrix algebra over a twisted group
algebra $k_\alpha(P\rtimes L)$ of the semidirect product $P\rtimes L$, for some
$\alpha\in H^2(L, k^\times)$, inflated to $P\rtimes L$. So there is a central $p'$-extension
\[ 1 \to Z \to H \to L \to 1 \]
and an idempotent $e$ in $kZ$, such that $k_\alpha(P\rtimes L)\cong$ $kGe$,
 where $G=P\rtimes H$. 

\begin{theorem}\label{th:qu}
With the notation and hypotheses above, let $\tilde\fA$ be a basic 
algebra of the twisted group algebra $k_{\alpha}(P \rtimes L)$. Then  
$k_{\alpha}(P \rtimes L)$ is a matrix algebra over $\tilde\fA$ and $\tilde\fA$ has 
an explicit presentation as a quantised version of the  group algebra 
$k(P\rtimes Z(H)/Z)$.
\end{theorem}

For the precise presentation and the proof, see Section \ref{se:ungrading}.
There are several new ingredients
required to extend the results from \cite{Benson/Kessar/Linckelmann:2019a}
to nonabelian defect groups. We first consider the associated graded 
$\gr_*(kGe)=\bigoplus_{n\geq 0}\ J^n(kGe)/J^{n+1}(kGe)$ of $kGe$,
briefly reviewed in the next section, and   make use of the 
Jennings--Quillen  theorem~\cite{Jennings:1941a,Quillen:1968a} and
Semmen \cite{Semmen:2005a}.
We show that $\gr_*(kGe)$   is 
isomorphic to a matrix algebra over 
a quantised version of the associated graded of the group algebra of the  
group  $P\rtimes (Z(H)/Z)$.  Specialising to  the  case $\alpha =0 $,  we  get a   presentation 
of  $\gr_*(k (P \rtimes L))$    which we have not seen before in the literature (see Remark~\ref{rk:fAquantised}).
The exact 
relations are stated  in Theorem~\ref{th:rels}; 
see also Theorems~\ref{th:Q/I} and~\ref{th:AxM} and
Corollary~\ref{co:basic}.
We then show that this  may be ungraded to exhibit the basic 
algebra of $kGe$ as a quantised version of the group algebra of
$P\rtimes (Z(H)/Z)$, see Section~\ref{se:ungrading}.

In Section~\ref{se:eg}, in order to illustrate the main results, 
we explicitly calculate the following examples of blocks with a normal extraspecial
defect group of order $p^3$ and exponent $p$ having a single isomorphism
class of simple modules. 

\begin{theorem} \label{pcubeexample}
Suppose that $p$ is odd.
Let $P$ be an extraspecial group of order $p^3$ and exponent $p$, let 
$H$ be a quaternion group of order $8$ acting on $P$  with $Z(H)$  acting trivially, and
with the two generators of $H$ inverting the two generators of $P$.
Set $G = P\rtimes H$.
The basic algebra of the associated graded $\gr_*(kP)$  of $kP$ is given by generators $x$, $y$, $z$,
subject to the relations
\[ x^p=0,\quad y^p=0,\quad xy-yx=z,\quad xz-zx=0,\quad yz-zy=0 \]
(these imply $z^p=0$)
while the basic algebra of the associated graded $\gr_*(kGe)$ of the non-principal block 
$e$ of $kG$ is 
given by generators $\x$, $\y$, $\z$, subject to the relations
\[ \x^p=0,\quad \y^p=0,\quad \x \y + \y \x = \z, \quad 
\x \z + \z \x = 0, \quad \y \z + \z \y = 0 \]
(these imply $\z^p=0$).
\end{theorem}

In the case $p=3$, we can be more precise and explicitly describe the algebra
$kP$ and a basic algebra of $kGe$ by `ungrading' the previous Theorem.

\begin{theorem}\label{27example}
With the notation of the previous theorem, assume that $p=3$.
The algebra $kP$ is given by generators $\tilde x$, $\tilde y$,
$\tilde z$, subject to the relations
\[ \tilde x^3=0,\quad \tilde y^3=0,\quad 
\tilde x\tilde y-\tilde y\tilde x=\tilde z,\quad 
\tilde x\tilde z-\tilde z\tilde x=\tilde z\tilde y\tilde z,
\quad \tilde y\tilde z-\tilde z\tilde y=-\tilde z\tilde x\tilde z \]
(these imply $\tilde z^p=0$)
while a basic algebra of  $kGe$ is given by generators $\tilde \x$, $\tilde \y$, $\tilde \z$,
subject to the relations
\[ \tilde\x^3=0, \quad \tilde\y^3=0,\quad 
\tilde\x\tilde\y+\tilde\y\tilde\x=\tilde\z, \quad
\tilde\x \tilde \z +
\tilde \z \tilde \x =
 -\tilde \z\tilde \y\tilde \z, \quad
\tilde \y \tilde \z + \tilde \z \tilde \y=
-\tilde \z \tilde \x \tilde \z \]
(these imply $\tilde\z^p=0$).
\end{theorem}

In Section  \ref{se:char2} we give an example in characteristic two, with
$P$ an extraspecial group of order $2^{1+4}$ and $H$ an
extraspecial group of order $3^{1+2}$.

Finally, the appendix contains some corrections 
to the calculations in 
\cite{Benson/Kessar/Linckelmann:2019a}.\bigskip

\noindent
{\bf Notation.}
The bracket $[-,-]$ is used in three different ways, depending on the context: 
as commutator $[g,h]=ghg^{-1}h^{-1}$ for elements $g$, $h$ in a 
multiplicatively written group, as Lie bracket in a Lie
algebra, and as additive commutator $[a,b]=ab-ba$ for elements $a$, $b$ in an
associative algebra.\bigskip

\noindent
{\bf Acknowledgements.} 
The first author is grateful to City, University of London for its
hospitality during the research for this paper, and 
to Ehud Meir for conversations
about the proof of Theorem~\ref{th:27}. The second author  acknowledges support from
EPSRC grant EP/T004592/1.

\section{The associated graded} 

The associated graded of a finite-dimensional $k$-algebra $A$ is the graded algebra
\[ \gr_*(A)=\bigoplus_{n\geq 0} J^n(A)/J^{n+1}(A), \] 
with the summands $J^n(A)/J^{n+1}(A)$
in degree $n$, where we adopt the convention 
$J^0(A)=A$. The image in $A/J(A)$ of a block idempotent of $A$ is a block idempotent
of $\gr_*(A)$, and this induces a bijection between the blocks of $A$ and the blocks 
of $\gr_*(A)$. Similarly, the image in $A/J(A)$ of a primitive idempotent of $A$ 
 is a primitive idempotent in $\gr_*(A)$. It follows  that $A$ and $\gr_*(A)$ have
 the same quiver.\medskip

Let $P$ be a finite $p$-group, $L$ an abelian $p'$-subgroup of $\Aut(P)$, and let
$\alpha\in H^2(L, k^\times)$. Since $k$ is algebraically closed, the canonical group
homomorphism $Z^2(G,k^\times)\to H^2(G,k^\times)$ splits (see for
example Theorem~11.15 of Isaacs~\cite{Isaacs:1976a}). Thus we may represent
$\alpha$ by a $2$-cocycle having the same order in $Z^2(G,k^\times)$ as its image in
$H^2(G,k^\times)$, still denoted by $\alpha$. Such a choice of $\alpha$ yields 
a central $p'$-extension
\[ 1 \to Z \to H \to L \to 1 \]
and a faithful character $\chi : Z \to k^\times$ such that $Z = [H,H]$ and such 
that for some choice of  inverse images $\hat x$ in $H$ for all $x$, we have 
$$\alpha(x,y)=\chi(\hat x\hat y \widehat{xy}^{-1})$$ 
for all $x$, $y\in L$. Moreover, $|Z|$ is equal to the order of $\alpha$ in $H^2(L,k^\times)$;
that is, the subgroup of $k^\times$ generated by the values of $\alpha$ is equal
to $\chi(Z)$.

Set $G = P \rtimes H$, where $H$ acts on $P$ via the canonical map $H\to L$, so that
$Z=$ $C_H(P)\leq Z(H)$, and hence $Z\leq Z(G)$. Thus the idempotent
\[ e=\frac{1}{|Z|}\sum_{z\in Z}\chi(z^{-1})z, \]
is a non-principal block of  $kG$, and the canonical surjection $G\to P\rtimes L$ 
with kernel $Z$  induces an algebra  isomorphism
$$kGe \xrightarrow{\cong} k_\alpha (P\rtimes L), $$
where $\alpha$ is inflated to $P\rtimes L$ via the canonical surjection 
$P\rtimes L\to L$.

We wish to describe $kGe$. This being difficult, we tackle first the associated 
graded algebra 
\[ \gr_*(kGe)=\bigoplus_{n\ge 0}J^n(kGe)/J^{n+1}(kGe). \]
Our goal is to give an explicit presentation of this as a 
quantum deformation of the corresponding associated graded 
for the (untwisted) group algebra $k(P \rtimes Z(H)/Z)$.

First, we recall the Jennings--Quillen theorem~\cite{Jennings:1941a,Quillen:1968a}
for the associated graded of $kP$. Our treatment follows Section~3.14 of~\cite{Benson:1998b}.
For $r\ge 1$, we have dimension subgroups 
\[ F_r(P)=\{g\in P \mid g-1 \in J^r(kP)\}. \]
Thus $F_1(P)=P$, $F_2(P)=\Phi(P)$, $[F_r(P),F_s(P)]\subseteq F_{r+s}(P)$, 
and if $g\in F_r(P)$ then $g^p\in F_{pr}(P)$. Furthermore, $F_r(P)$ is the
most rapidly descending central series with these properties. Define 
\[ \Jen_*(P) = \bigoplus_{r\ge 1}k\otimes_{\bF_p} F_r(P)/F_{r+1}(P). \]
Then $\Jen_*(P)$ is a $p$-restricted Lie algebra with Lie bracket induced by taking
 commutators in $P$  and
$p$th power map coming from taking $p$th powers  in $P$. As a restricted Lie algebra, 
$\Jen_*(P)$  is generated by its degree one elements because  the subgroups $F_r(P)$
 form the lowest central series with the properties mentioned above. 
Let $\cU\Jen_*(P)$ be
the restricted universal enveloping algebra of $\Jen_*(P)$ over $k$. 
As an associative algebra,
$\cU\Jen_*(P)$ is generated by its degree one elements. The commutator $[g,h]$ of
two elements $g$, $h\in P$ becomes the Lie bracket of the images of $g$, $h$ in 
$\Jen_*(P)$, and the image of that Lie bracket in $\cU\Jen_*(P)$ is in turn equal to 
the additive commutator  of the images of $g$, $h$ in the associative algebra 
$\cU\Jen_*(P)$.

The Jennings--Quillen theorem states that  there is  a  $k$-algebra  isomorphism
\[ \cU\Jen_*(P)\to \gr_*(kP) \] 
which  for any  $r$ and any   $g\in F_r(P)$  sends the image of $g $ in $F_r(P)/F_{r+1} (P) $ 
to the image of $g-1$ in $\gr_*(kP)$. 

The group action of $H$ on $P$ induces an action of $H$ on $\Jen_*(P)$ as
a restricted Lie algebra, because the Lie bracket in $\Jen_*(P)$ is induced by
taking commutators in $P$ and the $p$-power map in $\Jen_*(P)$ is induced by
taking $p$-th powers in $P$.  The Jennings--Quillen map is equivariant with respect to $H$, 
and therefore extends to an isomorphism
\[ \cU\Jen_*(P) \rtimes H  \xrightarrow{\ \cong\ } \gr_*(kP)\rtimes H 
\xrightarrow{\ \cong\ } \gr_*(kG), \]
where the second isomorphism uses the fact that $J(kG)=J(kP)kG=kGJ(kP)$ since
$H$ is a $p'$-group (cf. \cite[Theorem 4]{Semmen:2005a}). 
Since we have a canonical bijection between the blocks of $kG$ and the blocks
of $\gr_*(kG)$ as described at the beginning of this section, it follows that 
the blocks of  both $kG$ and $\gr_*(kG)$  are the idempotents in $kZ$.

\begin{remark}\label{rk:JenkGe}
If $e$ is an idempotent in $kH$, 
then the restriction of the projective module $kGe$ to $P$ is a direct sum
of $\dim_k(kHe)$ copies of $kP$. Furthermore, the radical layers of $kGe$ as a
$kG$-module are the same as the radical layers as a $kP$-module.
So we have
\begin{align*} 
\sum_{i\ge 0} \dim_k J^i(kGe)/J^{i+1}(kGe)&=\dim_k(kHe).
\sum_{i\ge 0}\dim_k J^i(kP)/J^{i+1}(kP) \\
&=\dim_k(kHe).\prod_r\left(\frac{1-t^{pr}}{1-t^r}\right)^{\dim_k\Jen_r(P)}. 
\end{align*}
It can also be seen by restriction to $P$ that if 
$e$ is a central idempotent in $kH$ 
then the associated graded $\gr_*(kGe)$ of the
algebra $kGe$ is generated by its degree zero and degree one elements.
\end{remark}

\begin{remark}
The algebra $\cU\Jen_*(P)$ is a 
finite dimensional cocommutative Hopf algebra, which defines a connected unipotent
finite group scheme $\cP$ whose group algebra is $k\cP\cong \cU\Jen_*(P)$. 
The finite group $H$ acts as automorphism on $\cP$, so we may form the semidirect product
$\cG = \cP \rtimes H$, which is again a finite group scheme.
\end{remark}

\section{The quantum relations} 

In this section, we define an algebra $\fA$, which will turn out to be a 
basic algebra for $\gr_*(kGe)$. The quantum commutation rules for $\fA$ are
given in Theorem~\ref{th:rels}, and the fact that $\fA$ is indeed a basic algebra
is shown in Corollary~\ref{co:basic}.

By \cite[Proposition 3.1]{Benson/Kessar/Linckelmann:2019a}
we have a bijection
\[ \Irr(Z(H)|\chi) \xrightarrow{\cong} \Irr(H|\chi),\qquad \phi
  \mapsto \tau_\phi \]
between one-dimensional characters of $Z(H)$ lying over $\chi$ and
irreducible characters of $H$ lying over $\chi$,
such that $\tau_\phi$ lies over $\phi$. The central idempotent corresponding
to $\tau_\phi$ is 
\[ e_\phi=\frac{1}{|Z(H)|}\sum_{h\in Z(H)}\phi(h^{-1})h. \]
Then $\displaystyle e = \sum_{\phi\in\Irr(Z(H)|\chi)} e_\phi$, and hence 
\[ kHe = \prod_{\phi\in\Irr(Z(H)|\chi)} kHe_\phi\ . \]
The factors $kHe_\phi$ are matrix algebras, corresponding to $\tau_\phi$,
all of the same dimension.  An element  $\xi $ of  $\Hom(H/Z, k^\times)$ induces 
an algebra automorphism  of $kHe$ sending
$he$ to $\xi(h)^{-1}he$.   This yields  an action of  $\Hom(H/Z, k^\times)$  on $kHe$ 
by   algebra automorphisms   which in   turn induces  a  permutation action    of $\Hom(H/Z, k^\times)$   on  the  set of 
  factors $kHe_\phi$.  The stabiliser of  any   factor is the subgroup  $\Irr(H/Z(H) ) $ of $\Irr(H/Z) $  and    elements  of 
  $\Irr(H/Z(H) ) $ act as inner automorphisms  on each factor.

Choose $\phi_0\in\Irr(Z(H)|\chi)$, and set $\tau=\tau_{\phi_0}$.
For each $\phi \in \Irr (Z(H)| \chi) $,
choose  a one dimensional representation  $\xi_\phi  \in  \Irr( H/Z)  $   inflated to $H$  whose
restriction to $Z(H)$ is $\phi\phi_0^{-1}$, and so 
that $\xi_{\phi_0}=1$.   The  $\xi_{\phi} $ form a set of coset 
representatives of $\Irr(H/Z(H) )$ in $\Irr( H/Z) $.
The algebra automorphism  induced by   $\xi_{\phi} $ sends $e_{\phi_0}$ to
$e_\phi$, hence restricts to an algebra isomorphism
$$kHe_{\phi_0} \cong kHe_\phi$$
sending $he_{\phi_0}$ to $\xi_\phi(h)^{-1}ee_\phi$. Taking the product over all $\phi$
yields a unital injective algebra homomorphism
$$kHe_{\phi_0} \to kHe$$
sending $he_{\phi_0}$ to $\sum_{\phi\in\Irr(Z(H)|\chi)}\ \xi_\phi(h)^{-1}he_\phi$.
By the above, this homomorphism depends on the choice of the $\xi_\phi$, but only up to
inner automorphisms of $kHe$.
We write $\fM$ for the image in $kHe$ of the matrix algebra $kHe_{\phi_0}$ under
this  algebra homomorphism. This is a unital matrix subalgebra in $kHe$.

We have a canonical homomorphism $\rho\colon H \to \Hom(H,k^\times)$
sending $g$ to $\rho(g)\colon h \mapsto \chi([h,g])$. The kernel of this homomorphism
is $Z(H)$ and its image is $\Hom(H/Z(H),k^\times)$. 
For each $\psi \in  \Irr (H/Z)  $    and each 
$\phi \in  \Irr (Z(H)|\chi) $, we write for simplicity $\phi\psi$ instead of
 $\phi\,(\psi\vert_{Z(H)})$. Then $\xi_{\phi\psi}\xi_\phi^{-1}\psi^{-1}$ is
trivial on $Z(H)$. So there exists an element $g_{\psi, \phi}\in H$ such that
\[ \rho(g_{\psi,\phi})=\xi_{\phi\psi}\xi_\phi^{-1}\psi^{-1}, \]
or equivalently, such that
$$\chi([h,g_{\psi,\phi}]) = \xi_{\phi\psi}(h)\xi_\phi(h)^{-1}\psi(h)^{-1}$$
for all $h\in H$.   We choose such elements $g_{\psi,\phi}$, one for each $\psi$ and $\phi$. Note that
these elements are unique up to multiplication by elements in $Z(H)$.

For any $\psi, \eta   \in     \Irr(H/Z)  $  and any   $\phi \in  \Irr (Z(H)| \chi) $,  we have 
$$ \rho (g_{\eta, \phi \psi}  g_{\psi, \phi} )=     \rho (g_{\eta, \phi \psi} )  \rho(g_{\psi, \phi}) =\xi_{\phi  \psi\eta }  \xi_{\phi}^{-1}  \eta^{-1} \psi^{-1}   =  \rho(g_{\eta\psi, \phi}).   $$

 \begin{lemma}\label{glemma}  
Let   $\psi_i,  \eta_j  \in  \Irr (H/Z) $, 
$\phi_i, \zeta_j \in \Irr(Z(H) |\chi) $, 
$ 1\leq i \leq m $, $ 1\leq j \leq n $.  Suppose that   
 $\phi_i  = \phi _{i-1} \psi_{i-1}  $ for all $ 2 \leq i \leq m $.   Then
 \begin{enumerate} 
 \item  $  g_{\psi_m, \phi_m} \ldots     
g_{\psi_1, \phi_1}   =     g _{\psi_m \cdots\psi_1,   \phi_1}   z $
 for some  $ z \in   Z(H)$. 
 \item  Suppose  further  that $\zeta_j = \zeta _{j-1} \psi_{j-1}  $ 
for all $ 2 \leq j\leq n $,   $\phi_1= \zeta_1  $  and   
$\psi_m \ldots  \psi_1 =   \eta_n \ldots  \eta _1 $. Then
 $  g_{\psi_m, \phi_m} \ldots     g_{\psi_1, \phi_1}   =   
g_{\eta_n, \zeta_n} \ldots     g_{\eta_1, \zeta_1}   z' $  for some  $z'  \in Z(H)$. 
 \end{enumerate}
 \end{lemma} 
\begin{proof}     
Since $Z =\Ker (\rho )$,    
(i)  follows by repeated application of the  equation 
displayed above the lemma.      Part  (ii)  follows  from (i)  applied 
to both   $  g_{\psi_m, \phi_m} \ldots     g_{\psi_1, \phi_1}   $ and 
$g_{\eta_n, \zeta_n} \ldots     g_{\eta_1, \zeta_1}  $.
\end{proof}

\begin{chunk}\label{chunk:PBW}

Since $k$ is algebraically closed, we may choose a $k$-basis $w_1,\dots,w_m$ of 
$\Jen_*(P)$, where $p^m=|P|$, consisting of homogeneous eigenvectors
of the action of $H$. We arrange the indices 
in such a way that if $i\le j$ then $\deg(w_i)\le \deg(w_j)$.
Then for each $w_i$ there is a character $\psi_i$ of $L$,
inflated to $H$, such that
\[ {^gw_i}=\psi_i(g)w_i \]
for $g\in H$. 
Define structure constants $c_{i,j,k}$ and $d_{i,k}$ for $\Jen_*(P)$ via
\[ [w_i,w_j]=\sum_k c_{i,j,k}w_k,\qquad w_i^{[p]}=\sum_k d_{i,k}w_k. \]
Here, $[w_i,w_j]$ denotes the Lie bracket and $w_i^{[p]}$ the
$p$-restriction map in $\Jen_*(P)$. 
We have 
\[ {^g[w_i,w_j]} = [{^gw_i}, {^gw_j}]  = \psi_i(g)\psi_j(g)[w_i,w_j], \]
\[ {^g(w_i^{[p]})} = (^gw_i)^{[p]} = (\psi_i(g)w_i)^{[p]} = \psi_i(g)^pw_i^{[p]}. \] 
It follows that 
if $c_{i,j,k}\ne 0$ then $\psi_i\psi_j=\psi_k$, and if $d_{i,k}\ne 0$ then 
$\psi_i^p=\psi_k$.

By the Poincar\'e--Birkhoff--Witt (PBW) theorem for
restricted Lie algebras (Jacobson~\cite{Jacobson:1962a}, page~190), the algebra 
$\cU\Jen_*(P)\cong \gr_*(kP)$ has a basis $\cB$
consisting of words $w_{i_1} \dots w_{i_r}$ where $i_1\le \dots\le
i_r$, and each index is repeated at most $p-1$ times
(so  we are  writing  $w_i^a $ as $w_i \dots w_i$). 
We follow the convention is that the empty word 
denotes the identity element in degree zero.
The element $w_{i_1}\dots w_{i_r}$ is an eigenvector 
for the conjugation action of
$H$, with character $\psi_{i_1}\dots\psi_{i_r}$.

In what follows we identify $\Jen_*(P)$ with its image in $\cU\Jen_*(P)\rtimes H$.
The calculations that follow are similar to those 
in Section~4 of~\cite{Benson/Kessar/Linckelmann:2019a} 
(with the corrections described  in Section~\ref{se:errata} below).   
For any  $\phi  \in \Irr(Z(H)|\chi) $, $w_i $ a  
basis element of $\Jen_*(P)$,   with  associated linear characters $\psi_i$, we write 
$g_{i, \phi}$  for the element  $g_{\psi_i,\phi} $.
\end{chunk}

\begin{lemma} \label{fALemma}
With the notation above,
the following equations in $(\cU\Jen_*(P)\rtimes H)e$ hold for all $h\in H$, all basis elements 
$w_i$ of $\Jen_*(P)$, the associated linear characters $\psi_i\in\Hom(H,k^\times)$  and 
all $\phi \in \Irr(Z(H)|\chi)$.

\begin{enumerate}
\item[{\rm (i)}]
$w_i e_\phi = e_{\phi\psi_i} w_i$.

\item[{\rm (ii)}]
$ (g_{i,\phi}w_i)(\xi_\phi(h)^{-1}e_\phi\, h) = (\xi_{\phi\psi_i}(h)^{-1}e_{\phi\psi_i}\, h)(g_{i,\phi}w_i). $

\item[{\rm (iii)}]
 $g_{i,\phi}w_ie_\phi=e_{\phi\psi_i}g_{i,\phi}w_i$ commutes with $\fM$.
\end{enumerate}
\end{lemma}

\begin{proof}
We have $hw_ih^{-1}=$ $\psi_i(h)w_i$, hence $w_i h= \psi_i(h)^{-1} hw_i$. Thus if $h\in Z(H)$, then
$\phi(h)^{-1}w_ih = \phi(h)^{-1}\psi_i(h)^{-1}hw_i$. Taking the sum over all $h\in Z(H)$ and
dividing by $|Z(H)|$ shows (i).  Note that $[g,h]e=$ $\chi([g,h])e$ for all $g$, $h\in H$.
Thus $g_{i,\phi}he=$ $\chi([h,g_{i,\phi}])^{-1} hg_{i,\phi}e$. It follows that
$$\xi_\phi(h)^{-1} g_{i\phi}w_ih e_\phi =\xi_\phi(h)^{-1} \psi_i(h)^{-1} g_{i,\phi}h w_i e_\phi=
\xi_\phi(h)^{-1}\psi_i(h)^{-1}\chi([h,g_{i,\phi}])^{-1} hg_{i,\phi} e_{\phi\psi_i} w_i\ ,$$
 where we have
used (i). Note that $e_{\phi\psi_i}$ is central in $kH$. Using the definition of $\rho$, the
scalar in the last expression is $\xi_{\phi\psi_i}(h)^{-1}$. This shows (ii). 
The equality in (iii) is the special case of (ii) applied with $h=1$. For the commutation with
$\fM$, we need to check that the elements in the statement commute with expressions of
the form $\sum_{\phi'} \xi_{\phi'}(h)^{-1}h e_{\phi'}$. This follows easily using (ii) and the fact
that the $e_\phi$ are pairwise orthogonal. 
\end{proof}

\begin{defn}\label{def:fA}
We define $\w_{i,\phi}=g_{i,\phi}w_ie_\phi$, and
let $\fA$ be  the subalgebra of 
$(\cU\Jen_*(P)\rtimes H)e$ generated by the elements $e_\phi$ and
$\w_{i,\phi}$. 
\end{defn} 

By Lemma \ref{fALemma}, the subalgebras $\fA$ and $\fM$ of $(\cU\Jen_*(P)\rtimes H)e$
commute.

\begin{lemma}\label{deg1genlemma}  The  algebra $\fA$ is generated by the  elements $e_\phi$ and
$\w_{i,\phi}$ for those  $i$ such that the element $w_i$ of $\Jen_*(P)$ has degree one.
\end{lemma} 

\begin{proof}  
 Since   $\Jen_*(P)$ is
generated  by  elements in degree one,   there exists a   basis  ${\mathcal  V}$  of $ \cU\Jen_*(P)$ consisting of    a subset of the set of  monomials  in   the   degree  one $w_i $'s.    Let $ w_t$  be  an arbitrary  element of the  chosen  basis  of  $\Jen_*(P) $
and write 
 $$ w_t =\sum_{ v\in  {\mathcal V}}    \alpha_v v.  $$  
If $u, u' \in  \cU\Jen_*(P)   $ are    eigenvectors  for the $H$ action corresponding to characters $\psi $ and   $\psi' $ respectively, then
$uu'$  is an  $H$-eigenvector  with   corresponding character $\psi \psi' $.
From this it follows that   if    a monomial   $   v= w_{i_m} \ldots
w_{i_1}   $  in degree  one  elements $w_{i_j} $    is an element of $
{\mathcal V}   $  such that  $\alpha _v \ne 0 $,  then  $ \psi_t =
\psi_{i_m}  \ldots   \psi_{i_1} $,   where for  each   $j$,   $1\leq j
\leq m$,  $\psi_{i_j}  \in\Irr( H/Z) $     is the  character   of $H$
corresponding to   the action on  $w_{i_j} $.     Let $\zeta  \in
\Irr( Z(H) |  \chi ) $  and let $v$ be as above.   By
Lemma~\ref{glemma},  
$$ g_{\psi, \zeta }   =    g_{i_m,  \phi_m}    \ldots  g_{i_1,  \phi_1 } z$$
where    $ z  \in Z(H) $,   $ \phi_1 =  \zeta $ and  $ \phi_j =   \phi_{j-1}\psi_{i_{j-1}}$, $ 2\leq j \leq m $.
On the other hand,    since   every $w_{i_j} $  is an  eigenvector for the $H$ action, 
$$ v  g_{i_m,  \phi_m}    \ldots  g_{i_1,  \phi_1 }     =\beta_v   w_{i_m}g_{i_m, \phi_m}   \ldots    w_{i_1}  g_{i_1, \phi_1}  $$  
for some $\beta_v\in k^{\times}  $.    Since  $z  e_{\zeta}  $  is a
non-zero scalar multiple of  $e_{\zeta} $, the above  equation   and
Lemma~\ref{fALemma} (iii) give that 
$$ v g_{\psi, \zeta}   e_{\zeta}  =  q_v \w_{i_m,  \phi_m}    \ldots  \w_{i_1,  \phi_1 } $$
for  some  non-zero scalar  $q_v$.  
Since all $w_{i_j}  $ are in degree one, it follows  that
$$ \w_{t, \zeta}    =   w  g_{\psi, \zeta}   e_{\phi} =  
\sum_{ v\in  {\mathcal V}}    \alpha_v v g_{\psi,  \zeta}
e_{\zeta}   $$ 
is  a linear combination of  monomials    in the  $\w_{i,  \phi} $  for those $i$ such that    $w_i $ has  degree one. 
\end{proof}

\begin{defn}
We define elements $z_{i,j,\phi}$, $z'_{i,j,k,\phi}$ and $z''_{i,k,\phi}$ in $Z(H)$
as follows.  By Lemma \ref{glemma}   we have 
\[ g_{j,\phi\psi_i}g_{i,\phi}=g_{i,\phi\psi_j}g_{j,\phi}z_{i,j,\phi} \]
for some $z_{i,j,\phi}\in Z(H)$.  If $c_{i,j,k}\ne 0$ then 
\[ g_{j,\phi\psi_i}g_{i,\phi}=g_{k,\phi}z'_{i,j,k,\phi} \]
for some $z'_{i,j,k,\phi}\in Z(H)$.
If $d_{i,k}\ne 0$ then 
\[ g_{i,\phi\psi_i^{p-1}}\dots 
g_{i,\phi\psi_i}g_{i,\phi}=g_{k,\phi}z''_{i,k,\phi} \] 
for some $z''_{i,k,\phi}\in Z(H)$.
\end{defn}

\begin{remark}
We have 
\[ z_{i,j,\phi}e_\phi=\phi(z_{i,j,\phi})e_\phi,\qquad 
z'_{i,j,k,\phi}e_\phi=\phi(z'_{i,j,k,\phi})e_\phi,\qquad
z''_{i,k,\phi}e_\phi=\phi(z''_{i,k,\phi})e_\phi. \] 
Also, we have $z_{i,j,\phi}z_{j,i,\phi}=1$, and
if
$c_{i,j,k}\ne 0$ then 
$z'_{i,j,k,\phi}=z'_{j,i,k,\phi}z_{i,j,\phi}$. 
\end{remark}

\begin{theorem}\label{th:rels}
Defining constants
\begin{align*} 
q_{i,j,\phi}&=\psi_i(g_{j,\phi}z_{i,j,\phi})\psi_j(g_{i,\phi}^{-1}z_{i,j,\phi})\phi(z_{i,j,\phi}), \\
q'_{i,j,k,\phi}&=\psi_j(g_{i,\phi})^{-1}\psi_k(z'_{i,j,k,\phi})\phi(z'_{i,j,k,\phi}), \\ 
q''_{i,k,\phi}&=\psi_i(g_{i,\phi\psi_i^{p-2}})^{-1}\dots\psi_i(g_{i,\phi\psi_i})^{-p+2}
\psi_i(g_{i,\phi})^{-p+1}\psi_k(z''_{i,k,\phi})\phi(z''_{i,k,\phi}) ,
\end{align*}
 we have
\begin{align} 
\w_{j,\phi\psi_i}\w_{i,\phi}-q_{i,j,\phi}
\w_{i,\phi\psi_j}\w_{j,\phi}&=
\sum_kc_{i,j,k}q'_{i,j,k,\phi}\w_{k,\phi} \label{eq:rels1} 
\\
\w_{i,\phi\psi_i^{p-1}}\dots \w_{i,\phi\psi_i}\w_{i,\phi}&=
\sum_kd_{i,k,\phi}\,q''_{i,k,\phi}\w_{k,\phi}. \label{eq:rels2}
\end{align}

By changing the choices of $g_{i,\phi}$ by elements of $Z(H)$, we may ensure                    
that $z_{i,j,\phi}\in Z$, and then the formula for the parameters $q_{i,j,\phi}$ simplifies to                     
\[ q_{i,j,\phi}=\psi_i(g_{j,\phi})\psi_j(g_{i,\phi}^{-1})\chi(z_{i,j,\phi}).  \]
\end{theorem}

\begin{proof}
We have
\begin{align*}
\w_{j,\phi\psi_i}\w_{i,\phi}&=
(g_{j,\phi\psi_i}w_je_{\phi\psi_i})(g_{i,\phi}w_ie_\phi)\\
&=g_{j,\phi\psi_i}w_jg_{i,\phi}w_ie_\phi \\
&=\psi_j(g_{i,\phi})^{-1}g_{j,\phi\psi_i}g_{i,\phi}w_jw_ie_\phi \\
&=\psi_j(g_{i,\phi})^{-1}g_{j,\phi\psi_i}g_{i,\phi}(w_iw_j+[w_i,w_j])e_\phi\\
&=\psi_j(g_{i,\phi})^{-1}g_{i,\phi\psi_j}g_{j,\phi}z_{i,j,\phi}w_iw_je_\phi \\
&\qquad{}+\psi_j(g_{i,\phi})^{-1}g_{j,\phi\psi_i}g_{i,\phi}[w_i,w_j]e_\phi\\
&=\psi_j(g_{i,\phi})^{-1}\psi_i(z_{i,j,\phi})\psi_j(z_{i,j,\phi})g_{i,\phi\psi_j}g_{j,\phi}w_iw_jz_{i,j,\phi}e_\phi \\
&\qquad{}+\sum_kc_{i,j,k}\psi_j(g_{i,\phi})^{-1}g_{k,\phi}z'_{i,j,k,\phi}w_ke_\phi \\
&=\psi_j(g_{i,\phi})^{-1}\psi_i(z_{i,j,\phi})\psi_j(z_{i,j,\phi})\psi_i(g_{j,\phi})\phi(z_{i,j,\phi})g_{i,\phi\psi_j}w_ig_{j,\phi}w_je_\phi \\
&\qquad{}+\sum_kc_{i,j,k}\psi_j(g_{i,\phi})^{-1}\psi_k(z'_{i,j,k,\phi})g_{k,\phi}w_kz'_{i,j,k,\phi}e_\phi\\
&=\psi_i(g_{j,\phi}z_{i,j,\phi})\psi_j(g^{-1}_{i,\phi}z_{i,j,\phi})\phi(z_{i,j,\phi})(g_{i,\phi\psi_j}w_ie_{\phi\psi_j})(g_{j,\phi}w_je_\phi) \\
&\qquad{}+\sum_kc_{i,j,k}\psi_j(g_{i,\phi})^{-1}\psi_k(z'_{i,j,k,\phi})\phi(z'_{i,j,k,\phi})g_{k,\phi} w_ke_\phi\\
&=q_{i,j,\phi}\w_{i,\phi\psi_j}\w_{j,\phi}
+\sum_kc_{i,j,k}q'_{i,j,k,\phi}\w_{k,\phi}.
\end{align*}

Similarly,
\begin{align*}
\w_{i,\phi\psi_i^{p-1}}&\dots \w_{i,\phi\psi_i}\w_{i,\phi} \\
&=
(g_{i\phi\psi_i^{p-1}}w_ie_{\phi\psi_i^{p-1}})\dots
(g_{i,\phi\psi_i}w_ie_{\phi\psi_i})(g_{i,\phi}w_ie_\phi)\\
&=g_{i,\phi\psi_i^{p-1}}w_i\dots g_{i,\phi\psi_i}w_ig_{i,\phi}w_ie_\phi \\
&=\psi_i(g_{i,\phi\psi_i^{p-2}})^{-1}\dots\psi_i(g_{i,\phi\psi_i})^{-p+2}
\psi_i(g_{i,\phi})^{-p+1}(g_{i,\phi\psi_i^{p-1}}\dots g_{i,\phi\psi_i}g_{i,\phi})w_i^pe_\phi\\
&=\sum_kd_{i,k}\psi_i(g_{i,\phi\psi_i^{p-2}})^{-1}\dots\psi_i(g_{i,\phi\psi_i})^{-p+2}
\psi_i(g_{i,\phi})^{-p+1}g_{k,\phi}z''_{i,k,\phi}w_ke_\phi\\
&=\sum_kd_{i,k}\psi_i(g_{i,\phi\psi_i^{p-2}})^{-1}\dots\psi_i(g_{i,\phi\psi_i})^{-p+2}
\psi_i(g_{i,\phi})^{-p+1}\psi_k(z''_{i,k,\phi})g_{k,\phi}w_kz''_{i,k,\phi}e_\phi\\
&=\sum_kd_{i,k}\psi_i(g_{i,\phi\psi_i^{p-2}})^{-1}\dots\psi_i(g_{i,\phi\psi_i})^{-p+2}
\psi_i(g_{i,\phi})^{-p+1}\psi_k(z''_{i,k,\phi})\phi(z''_{i,k,\phi})g_{k,\phi}w_ke_\phi\\
&=\sum_k d_{i,k}q''_{i,k,\phi}\w_{k,\phi}.
\end{align*}
For the final remark, 
just as in Lemma~4.12\,(3) of~\cite{Benson/Kessar/Linckelmann:2019a}, 
we may change the choices of $g_{i,\phi}$ by elements of $Z(H)$ to 
ensure that $z_{i,j,\phi}\in Z$, with the same argument. Then the characters $\psi_i$ take value one 
on these elements, leading to the given simplifications of the constants.
\end{proof}

Recall   that by Lemma~\ref{deg1genlemma}, $ \fA$   is generated by the $e_{\phi}$  and    the  $\w_{i, \phi} $  for those $i$  such that     the element  $w_i$  of $\Jen_*(P)$ has degree one. 
\begin{theorem}\label{th:Q/I}
The algebra $\fA$ is given as a quiver with relations $kQ/I$, 
where $Q$ is 
the quiver with $|Z(H):Z|$ vertices labelled $[\phi]$ corresponding to the 
idempotents $e_\phi \in kZ(H)$ lying over $\chi$ and directed edges 
\[ [\phi] \xrightarrow{\quad i\quad }[\phi\psi_i] \]
corresponding to the  element 
\[ \w_{i,\phi}=g_{i,\phi}w_ie_\phi = e_{\phi\psi_i}g_{i,\phi}w_i=
e_{\phi\psi_i}g_{i,\phi}w_ie_\phi \]
for those $i$ such that the element $w_i$ of $\Jen_*(P)$ has degree one. The 
relations are those that follow from the structure constant relations of 
Theorem~\ref{th:rels},  where  for  each    $k$   such that  $w_k $   is in degree  greater than or equal to $2$,    any  $\w_{k, \zeta} $  appearing in Theorem~\ref{th:rels}   is replaced    by an element      in $kQ$ corresponding via Lemma~\ref{deg1genlemma}  to  an expression  for $\w_{k, \zeta}   $  in terms of the $\w_{i, \phi}  $ such that  $ w_i $ has  degree one. There is a PBW style basis $\cB'$ for $\fA$ (described below in
the proof), consisting of composable
monomials in the $\w_{i,\phi}$, giving $\dim(kQ/I)=\dim(\fA)=|Z(H):Z|\cdot|P|$.
\end{theorem}
\begin{proof}
By  Lemma~\ref{deg1genlemma}, $\fA$ is generated by the idempotents
$e_\phi$ and the elements $\w_{i,\phi}$. By
Lemma~\ref{fALemma} and Theorem~\ref{th:rels} they satisfy the
given relations. 
Thus we have a surjective homomorphism from $kQ/I$ to $\fA$ taking $[\phi]$
to $e_\phi$ and $[\phi] \xrightarrow{\quad i \quad}[\phi\psi_i]$ to
$\w_{i,\phi}$.

The relations holding in $kQ/I$ allow us to write every element of $\fA$ as a
linear combination of elements of the set $\cB'$ consisting of 
the $e_\phi$ and composable monomials
in the $\w_{i,\phi}$ where the indices $i$ are in order, and
each index $i$ is repeated at most $p-1$ times. The number of such
monomials (including the $e_\phi$) is $|Z(H):Z|\cdot  |P|$.    Replacing   $[\phi] $  by $e_\phi $, 
$\w_{i, \phi} $  by $[\phi] \xrightarrow{\quad i \quad}[\phi\psi_i]$   for those   $i $ such that $ w_i$ has degree one  and    $\w_{i, \phi} $  by their  chosen  lifts in $kQ$  for those   $i $ such that $ w_i$  has degree greater than or equal to  two, we  see   by the same  reasoning that  $\dim_kQ/I  $  is at most  $|Z(H):Z|\cdot  |P|$.  

If there were a linear relation in $\fA$ between the monomials in $\cB'$,
then there would be a linear relation between the ones of maximal
length, namely length $m(p-1)$. 
There is one of these for each $\phi$, and they are linearly
independent elements of the socle of $kG$ because they are non-zero
elements of different projective summands $kGe_\phi$. Thus $\dim(\fA)$
is equal to $|Z(H):Z|\cdot |P|$ and $kQ/I \to \fA$ is an isomorphism.
\end{proof}

\begin{rk}\label{rk:fAquantised}
The group algebra of the semidirect product 
$P\rtimes Z(H)/Z$, with the action given by
restricting the action of $H/Z$ on $P$, has only one block. We can
perform the computations above for this group, and the results look
similar, except that the factors of
$\phi(z_{i,j,\phi})$, $\phi(z'_{i,j,k,\phi})$, and $\phi(z''_{i,k,\phi})$ 
in the definitions of $q_{i,j,\phi}$, $q'_{i,j,k,\phi}$, and
$q''_{i,k,\phi}$ are missing in Theorem~\ref{th:Q/I}. So removing these factors, the
relations in Theorem~\ref{th:rels} are the relations in 
$\gr_*(k(P\rtimes Z(H)/Z))\cong \cU\Jen_*(P)\rtimes Z(H)/Z$.
Thus we can see $\fA$ as a quantum deformation
of the algebra $\cU\Jen_*(P)\rtimes Z(H)/Z$.  Also, we note that   in the case that  $\alpha=0$,  we have  $Z=1 $, $ H=L$   and  Theorem~\ref{th:Q/I}  provides an explicit presentation of  $\gr_*( k(P\rtimes L ))  $. 
\end{rk}

As in~\cite{Benson/Kessar/Linckelmann:2019a}, we now make use of the following lemma
(see Chapter~3, Corollary~4.3 in Bass~\cite{Bass:1967a}).

\begin{lemma}\label{le:Bass}
Let $A\le B$ be $k$-algebras with $A$ an Azumaya algebra (that is, a finite-dimensional
central separable $k$-algebra). Then the map $A \otimes_k C_B(A) \to B$ is an
isomorphism.\qed
\end{lemma}

\begin{theorem}\label{th:AxM}
The multiplication in $(\cU\Jen_*(P)\rtimes H)e$ induces an isomorphism 
$$\fA\otimes_k \fM \xrightarrow{\ \cong\ } (\cU\Jen_*(P)\rtimes H)e .$$
\end{theorem}

\begin{proof}
The proof is similar to that of  Theorem~4.15 of~\cite{Benson/Kessar/Linckelmann:2019a}.
Applying Lemma~\ref{le:Bass} with $A=\fM$ and $B$ the subalgebra generated by
$\fA$ and $\fM$, we see that the given map is injective. The dimensions are given
by $\dim(\fA)=|Z(H):Z|\cdot |P|$, $\dim(\fM)=|H:Z(H)|$ and 
$\dim((\cU\Jen_*(P)\rtimes H) e)=$ $\dim(kGe)=|G:Z|$,
so $\dim((\cU\Jen_*(P)\rtimes H)e)=\dim(\fA)\cdot\dim(\fM)$ and the map is an isomorphism.
\end{proof}

\begin{cor}\label{co:basic}
We have 
\[ \gr_*(kGe)\cong (\cU\Jen_*(P)\rtimes H)e \cong \Mat_m(\fA), \] 
where $m=\sqrt{|H:Z(H)|}$. In particular, $\fA$ is a  basic algebra of $\gr_*(kGe)$.\qed
\end{cor}

\begin{cor}
The algebra $\fA$ is generated by its degree zero and degree one
elements.
\end{cor}
\begin{proof}
This follows from Corollary~\ref{co:basic} and Remark~\ref{rk:JenkGe}.
\end{proof}

\section{Ungrading the relations}\label{se:ungrading} 

We saw in the last section that the relations for the basic algebra of
$\gr_*(kGe)$ are a quantised version of the relations for
$\gr_*(kP\rtimes(Z(H)/Z))$. In this section, we show that the same
holds without taking the associated graded.

Since $|H|$ is coprime to $p$, the characteristic of $k$, we can
choose invariant complements to $J^{n+1}(kP)$ in $J^n(kP)$ for each
$n\ge 0$. 
Let $w_1,\dots,w_m$ be  the  basis of   $\Jen_*(P)$  
chosen in Section~\ref{chunk:PBW}, and let $\cB$ be the resulting PBW
basis of $\cU\Jen_*(P)\cong \gr_*(kP)$ described there.
Regarding $\Jen_*(P)$ as an $k$-linear subspace of
$\gr_*(kP)$,
this enables us to choose representatives $\tilde w_i$ in $kP$ of
the $w_i$ in such a way that 
\[ g\tilde w_i g^{-1} = \psi_i(g)\tilde w_i. \]    
Let $\tilde \cB$ be the  corresponding  basis  of 
$kP$ consisting of monomials  in  the  $\tilde w_i $.  That  is,   
if $w_{i_1}   \dots  w_{i_r}  $  is  an element of $\cB$   
then the corresponding  element of $\tilde \cB $ is 
$\tilde w_{i_1}  \dots  \tilde w_{i_r} $.  
An element   $\tilde w_{i_1}  \ldots   \tilde w_{i_r} $    
of $\tilde\cB$ is an   eigenvector   for   the  action of  $H$  for the 
character  $\psi_{i_1}  \ldots  \psi_{i_r}  $.

When we ungrade a relation of the form
$[w_i,w_j]=\sum_kc_{i,j,k}w_k$,
we obtain a relation of the form
\begin{equation}\label{eq:[xi,xj]} 
[\tilde w_i,\tilde w_j] = \sum_kc_{i,j,k}\tilde w_k + y_{i,j} 
\end{equation}
in $kP$, where $y_{i,j}$ is a linear combination of elements of 
$\tilde\cB$ in a higher power of  the radical than $\deg(w_i)+\deg(w_j)$.  
Moreover,  each basis  monomial  
$\tilde w_{i_1}   \ldots \tilde w_{i_r} $   that  occurs in
$y_{ij}$ is an eigenvector  for the character   $\psi_i \psi_j $,  and  consequently
\begin{equation} \label{eq:ijproduct}  
\psi_{i_1}   \ldots  \psi_{i_r}   =\psi_i \psi_j .
\end{equation}
 Similarly, when we ungrade a relation of the form 
$w_i^p=\sum_kd_{i,k}w_k$, we obtain a relation of the form
\begin{equation}\label{eq:xi^p} 
\tilde w_i^p=\sum_kd_{i,k}\tilde w_k + y''_{i} 
\end{equation}
in $kP$, where $y''_{i}$ is    a linear  combination    of   monomial  basis elements  
in a higher power of the radical than $p.\deg(w_i)$.  
Each basis  monomial  $\tilde w_{i_1}\dots\tilde w_{i_r} $   that  occurs in 
 $y''_{i} $ is an eigenvector  for the character   $\psi_i ^p$  and  consequently
 \begin{equation} \label{eq:xi^pproduct}   
\psi_{i_r}  \ldots  \psi_{i_1} =\psi_i ^p. 
\end{equation}

\begin{defn}
As in Definition~\ref{def:fA}, we define 
$\tilde \w_{i,\phi}=g_{i,\phi}\tilde w_ie_\phi$. Then $\tilde \w_{i,\phi}$
commutes with $\fM$. We define $\tilde\fA$ to be the
subalgebra of $kGe$ generated by the elements  $e_\phi $  and  $\tilde \w_{i,\phi}$.  
For  an  element  $\tilde w=\tilde w_{i_1}   \ldots   \tilde w_{i_r} $
of   $\tilde\cB$  and a character $\phi\in\Irr(Z(H)|\chi)$, set 
$\tilde \w_{\phi}  = \tilde \w_{i_1,    \phi \psi_1\ldots \psi_r  }\dots
\tilde \w_{  i_{r-1},  \phi \psi_r} \tilde  \w_{i_r, \phi}  $. Denote by
$\tilde \cB'$ the subset of $\tilde\fA$ consisting of the elements
$\tilde x_\phi$ for $\tilde w$ in $\tilde\cB$ and
$\phi\in\Irr(Z(H)|\chi)$. We shall see below in Theorem~\ref{th:grfA} that
$\tilde\cB'$ is a basis for $\tilde\fA$.
\end{defn}

\begin{prop} \label{p:ungradedrelations}
The elements $\tilde \w_{i,\phi}$ satisfy the following relations.
\begin{align*} 
\tilde \w_{j,\phi\psi_i}\tilde \w_{i,\phi}-q_{i,j,\phi}\tilde
\w_{i,\phi}\tilde \w_{j,\phi\psi_i} &= \sum_kc_{i,j,k}q'_{i,j,k,\phi}\tilde
\w_{k,\phi} +\psi_j(g_{i,\phi})^{-1}g_{j,\phi\psi_i}g_{i,\phi}y_{i,j}e_\phi. \\
\tilde \w_{i,\phi\psi_i^{p-1}}\dots \tilde \w_{i,\phi\psi_i}\tilde
  \w_{i,\phi} &= \sum_k d_{i,k}q''_{i,k,\phi}\tilde \w_{k,\phi} 
+
\psi_i(g_{i,\phi\psi_i^{p-2}})^{-1}\dots\\
&\qquad\qquad\dots\psi_i(g_{i,\phi\psi_i})^{-p+2}\psi_i(g_{i,\phi})^{-p+1} 
(g_{i,\phi\psi_i^{p-1}}\dots g_{i,\phi\psi_i}g_{i,\phi})y''_ie_\phi.
\end{align*}
 Moreover, suppose that  $y_{i,j}  = \sum_{\tilde w \in \tilde\cB}
 c_{i,j,\tilde w}  \tilde w$ 
and     $ y_i''= \sum _{\tilde w \in \tilde\cB}   d_{i,\tilde w}
\tilde w $. 
For each $\tilde w \in \tilde\cB$,  there  exist    elements
$q'_{i,j,\tilde w, \phi}   $ 
and    $ q''_{i,\tilde w,\phi} $  of  $ k ^{\times} $ such that 
 \begin{align*}  g_{j,\phi\psi_i}g_{i,\phi}y_{i,j}e_\phi=   
\sum_{\tilde w \in \tilde\cB}  q'_{i,j,\tilde w,\phi} c_{i,j,\tilde w}  \tilde \w_\phi.  \\
 (g_{i,\phi\psi_i^{p-1}}\dots g_{i,\phi\psi_i}g_{i,\phi})y''_ie_\phi
   =    
\sum_{\tilde w \in \tilde\cB}  q''_{i,\tilde w,\phi} d_{i,\tilde w}  \tilde \w_\phi. 
  \end{align*}
  \end{prop}
\begin{proof}
Following through the proof of relation~\eqref{eq:rels1},
we can replace each $w$ with $\w$ until the sixth line,
where we have to use~\eqref{eq:[xi,xj]} for the commutator.
At this point, the extra term is 
\[ \psi_j(g_{i,\phi})^{-1}g_{j,\phi\psi_i}g_{i,\phi}y_{i,j}e_\phi. \]

Similarly, following through the proof of relation~\eqref{eq:rels2},
we can replace each $w$ with $\w$ until the fourth line,
where we have to use~\eqref{eq:xi^p}. At this point, the extra
term is
\begin{equation*}
\psi_i(g_{i,\phi\psi_i^{p-2}})^{-1}\dots
\psi_i(g_{i,\phi\psi_i})^{-p+2}\psi_i(g_{i,\phi})^{-p+1}
(g_{i,\phi\psi_i^{p-1}}\dots g_{i,\phi\psi_i}g_{i,\phi})y''_ie_\phi.
\end{equation*}  
By Lemma \ref{glemma}  and  Equation \ref{eq:ijproduct},  
for each $\w =  \w_{i_r}   \ldots   \w_{i_1} \in \tilde\cB$ such that
$c_{i,j, \tilde w, \phi} \ne 0 $,    
\[ g_{j,\phi\psi_i}g_{i,\phi}=  g_{i_r, \phi \psi_1\ldots \psi_r  }  
\ldots      g_{  i_2,  \phi \psi_1}   g_{i_1, \phi} z \]   
for some $z \in  Z(H) $.   The second assertion follows   
from this  by the   fact  that for any $g\in H$, any $z \in Z(H)$,  
any $\tilde w_i$, and any $\zeta \in \Irr( Z(H)  |  \chi) $,   
$g\tilde w_i   =\psi_i(g)\tilde w_i   g $  is a scalar multiple of
$\tilde w_i g$,  
$ ze_{\zeta} =\zeta(z) e_\zeta $ is a scalar multiple of 
$e_\zeta $   and   $g\tilde w_i e_{\zeta}  =  e_{\zeta\psi_i}  g
\tilde w_i e_\phi $. The last assertion  follows  in a similar fashion   from  
Lemma \ref{glemma}  and  Equation \ref{eq:xi^pproduct}.
\end{proof}

\begin{theorem}\label{th:grfA}
The algebra $\tilde\fA$ is given as a quiver with relations $kQ/\tilde
I$, where $Q$ is as in Theorem~\ref{th:Q/I}, but with edges
corresponding to the lifts $\tilde \w_{i,\phi}$ of the $\w_{i,\phi}$
given there. The relations are those that follow from the structure
constant relations of Proposition~\ref{p:ungradedrelations}, together
with relations saying that every composite of at least $s$ arrows is
zero, where $s$ is the radical length of $kP$.

The set $\tilde\cB'$ is a
PBW style basis of $\tilde\fA$, giving
$\dim\tilde\fA= |Z(H):H|$.
There is a natural isomorphism $\gr_*(\tilde \fA)\cong \fA$, sending 
each $\tilde \w_{i,\phi}$ to $\w_{i,\phi}$.
\end{theorem}
\begin{proof}
It follows from the relations in Proposition~\ref{p:ungradedrelations}
that the linear span of $\tilde\cB'$ is closed under multiplication
modulo a large enough power of the arrow ideal. The zero relations
for composites of $s$ arrows then show that this ideal is zero,
and therefore that $\tilde\cB'$ linearly spans $\tilde\fA$. 
The image of an element $\tilde \w_{i,\phi}$ in $\gr_*(kGe)$ is
equal to $\w_{i,\phi}$, which lies in $\fA$. Since the elements
$\w_{i,\phi}$ of $\cB'$ are linearly independent, it follows
that the elements $\tilde \w_{i,\phi}$ of $\tilde\cB'$ are linearly
independent, and therefore form a basis for $\tilde\fA$. This
therefore induces a natural isomorphism
$\gr_*(\tilde\fA)\cong \fA$. Since $\fA$ is generated by its degree
one elements, $\tilde\fA$ has the same quiver, with the lifts of the
relations. 
\end{proof}

\begin{theorem}\label{th:tildeAisbasic}
The multiplication in $kGe$ induces an isomorphism $\tilde\fA\otimes_k
\fM \to kGe$.
\end{theorem}
\begin{proof}
By Theorem~\ref{th:grfA} we have 
$\dim(\tilde\fA)=|Z(H):Z|\cdot|P|$. So this is now proved in the same way as
Theorem~\ref{th:AxM}. 
\end{proof}

\begin{cor}
We have $kGe\cong \Mat_m(\tilde\fA)$, where $m=\sqrt{|H:Z(H)|}$, so that
$\tilde \fA$ is the basic algebra of $kGe$.\qed
\end{cor}

\begin{rk}
As in Remark~\ref{rk:fAquantised}, if we perform the computations of
this section with the group algebra of the semidirect product
$P\rtimes Z(H)/Z$ instead of $kGe$, the results look similar except
with different scalars. So we can see $\tilde\fA$ as a quantum
deformation of the algebra $k(P\rtimes Z(H)/Z)$. This observation,
together with Theorems~\ref{th:grfA} and~\ref{th:tildeAisbasic},
complete the proof of Theorem~\ref{th:qu}.
\end{rk}

We shall see some explicit examples of the ungrading of the  
relations  in Section~\ref{se:eg}.

\section{\texorpdfstring{Example: $P$ extraspecial of order $p^3$ and
    exponent $p$}{Example: P extraspecial of order p³ and exponent p}}
\label{se:eg}

Let $k$ have characteristic $p$, an odd prime, and $P$ be an extraspecial
$p$-group of order $p^3$ and exponent $p$, with presentation
\[ P = \langle g,h,c \mid g^p=h^p=c^p=1,\ [g,h]=c,\ [g,c]=[h,c]=1\rangle. \]
We denote by $H$ the quaternion group of order $8$, given by a presentation
\[ H = \langle s, t\ |\ s^4=1,\ s^2=t^2,\ ts = s^{-1}t\rangle\cong Q_8. \]
Set $Z= \langle s^2\rangle$; this is the centre of $H$. 
We consider the following action of $H$ on $P$, and set $G=P\rtimes H$.
\[ g^s = g^{-1},\qquad g^t = g,\qquad 
h^s = h,\qquad h^t  = h^{-1} \]
It follows that $c^s=c^t=c^{-1}$, and $Z$ acts trivially on $P$.  
This action lifts the action of $C_2\times C_2$
on $C_p\times C_p\cong P/\langle c\rangle$, where here the nontrivial element of
each copy of $C_2$ acts as
inversion on the corresponding copy $C_p$. 
The group algebra $kG$ has two blocks,
namely the principal block $e_0 = \frac{1}{2} (1+s^2)$ and the nonprincipal block
$e = \frac{1}{2} (1-s^2)$ corresponding to the faithful central character 
$\chi\colon Z \to k^\times$ given by $\chi(s^2)=-1$. We shall be interested in $kGe$. 

\begin{rk}
Let $x=g-1$, $y=h-1$, $z=c-1$ in $kP$. Then
\begin{equation}\label{eq:Z} 
z = (xy-yx)(1+x)^{-1}(1+y)^{-1} 
\end{equation}
and a presentation for $kP$ is given by generators $x$ and $y$, and relations
saying that $x^p=0$, $y^p=0$, and the element $z$ defined by~\eqref{eq:Z} is central
with $p$th power equal to zero. 
Note that the element $(1+x)^{-1}(1+y)^{-1}$ is congruent
to $1$ modulo $J(kP)$, and so in the associated graded $\gr_*(kP)$ this
term in~\eqref{eq:Z} may be ignored. This is used in the proof of
Theorem \ref{pcubeexample} that follows.
\end{rk}

\begin{proof}[Proof of Theorem \ref{pcubeexample}]
Denote by $x$, $y$, $z$ the images of $g$, $h$, $c$ in $\Jen_*(P)$,
respectively. (These elements are mapped to the images of $g-1$, $h-1$, $c-1$ in
$\gr_*(kP)$ under the canonical map $\Jen_*(P)\to \gr_*(kP)$).
The three dimensional $p$-restricted Lie algebra  $\Jen_*(P)$ is spanned
by  the elements $x$, $y$ in degree one together with $z=[x,y]$ (by the previous
Remark) in degree two, satisfying $[x,y]=[y,z]=0$. The
$p$-restriction map given by $x^{[p]}=y^{[p]}=z^{[p]}=0$. 
Its $p^3$ dimensional universal enveloping algebra $\cU\Jen_*(P)$ is isomorphic to
$\gr_*(kP)$. This shows the first part of Theorem \ref{pcubeexample}.

The action of $H$ on $\Jen_*(P)$ is given by 
\[ x^s= - x,\qquad x^t = x,\qquad\qquad
y^s= y,\qquad y^t = -y,\qquad\qquad
z^s = -z,\qquad z^t = -z.  \]
The elements $x$, $y$ and $z$ are eigenvectors for $H$
on $\Jen_*(P)$. So we set 
$w_1=x$, $w_2=y$, $w_3=z$. The characters $\psi_i$ of $H$ satisfying
$gw_ig^{-1}=\psi_i(g)w_i$ for $g\in H$ are given as follows.
\[ \psi_1(s)=-1,\quad \psi_1(t)=1,\qquad
\psi_2(s)=1,\quad \psi_2(t)=-1,\qquad
\psi_3(s)=-1,\quad \psi_3(t)=-1. \]
Note that the relation $[x,y]=z$ in $\Jen_*(P)$ implies $\psi_1\psi_2=\psi_3$. 

Denoting as above by $e = \frac{1}{2} (1-s^2)$ the nonprincipal block of $kG$,
the block algebra  $kGe$ has a unique isomorphism class of simple modules.
Indeed, $e$ corresponds to the unique $2$-dimensional simple $kH$-module, and
hence the semisimple quotient of $kGe$ is the matrix
algebra $\fM=kHe\cong \Mat_2(k)$.

Since $Z=Z(H)$, there is only one central character of $Z(H)$ lying above $\chi$,
namely $\phi=\chi$, and $\xi_\phi=1$. The map 
$\rho\colon H/Z(H) \to \Hom(H/Z(H),k^\times)$ takes $s$ to $\phi_2$,
$t$ to $\phi_1$ and $st$ to $\phi_3$. 
Thus $g_{1,\phi}=t$, $g_{2,\phi}=s$ and $g_{3,\phi}=st$; these are only well
defined up to multiplication by $Z(H)$.

The block algebra $\gr_*(kGe)$ of $\gr_*(kG)$ also has one isomorphism class of simple 
modules, namely the same $2$-dimensional simple $kH$-module as above, and by Theorem~\ref{th:AxM} and Corollary~\ref{co:basic} we have 
\[ \gr_*(kGe) \cong  
\fA \otimes_k \fM\cong \Mat_2(\fA), \] 
where $\fM=kHe\cong \Mat_2(k)$ and $\fA=(\gr_*(kGe))^H$.
The algebra
$\fA$ contains elements $g_{1,\phi}w_1e=txe$, $g_{2,\phi}w_2e=sye$ and
$g_{3,\phi}w_3e=stze$. The constants are given by $q_{1,2,\phi}=-1$ and $q'_{1,2,3,\phi}=-1$, 
so these satisfy the following relation:
\[ (txe)(sye)+(sye)(txe) = -tsxye - styxe = st(xy-yx)e=stze \]
Similar computations give
\[ (txe)(stze) + (stze)(txe) = 0, \qquad (sye)(stze) + (stze)(sye) = 0. \]
Writing $\x = txe$, $\y = sye$ and $\z = stze$, we therefore have
\[ \x \y + \y \x = \z, \quad 
\x \z + \z \x = 0, \quad \y \z + \z \y = 0,\quad
\x^p=0,\quad \y^p = 0, \quad \z^p =0. \]
This is a presentation for the basic algebra $\fA$ of $\gr_*(kGe)$, with generators 
$\x$ and $\y$,
and with $\z$ defined as $\x \y + \y \x$.
This proves Theorem \ref{pcubeexample}.
\end{proof}

\begin{rk}
The first part of the above proof shows that $\Jen_*(P)$
is isomorphic to the $p$-restricted Lie algebra of $3\times 3$ matrices of the
form $\left(\begin{smallmatrix} 0 & * & * \\ 0 & 0 & * \\ 0 & 0 & 0 \end{smallmatrix}\right)$.
\end{rk}

In order to prove Theorem \ref{27example},
ungrading the algebra is our next task. The problem is that the
generators $g-1$ and $h-1$ of $kP$ are not well suited to
dealing with automorphisms. We have an
action of $\bF_p^\times \times \bF_p^\times$ on $P$ where $(i,j)$
sends $g$ to $g^i$ and $h$ to $h^j$. The commutator $c=[g,h]$
is sent to $c^{ij}$.  Set  
 \[ \tilde x=-\sum_{i=1}^{p-1}g^i/i,\qquad \tilde y=-\sum_{j=1}^{p-1}h^j/j. \] 

\begin{lemma}  \label{liftlemma}  We have  $\tilde x  \equiv    g-1  \pmod {J^2(kP)} $    and   $\tilde y  \equiv    h-1 \pmod {J^2(kP)}$.
\end{lemma}  

\begin{proof}   Since $ p$ is odd, $\sum_{i=1}^{p-1} 1/i =  \sum_{i=1}^{p-1} i= 0  $  in $k$  whence   $\tilde x  =  -\sum_{i=1}^{p-1}(g^i-1)/i $.
Now the assertion  for  $\tilde x $   follows   since   $(g^i-1)/i   \equiv     g-1   \pmod {J^2 (kP) }   $  for any $i$, $ 1\leq  i \leq p-1 $. The proof   for $\tilde y $ is similar.
\end{proof}  

Note that  $\tilde x$ an eigenvector in the $(1,0)$ eigenspace and
$\tilde y$ an eigenvector in the $(0,1)$ eigenspace 
of $\bF_p^\times \times \bF_p^\times$.
Then we set $\tilde z=[\tilde x,\tilde y]=\tilde x\tilde y-\tilde
y\tilde x$, an eigenvector in the
$(1,1)$ eigenspace.   By Lemma~\ref{liftlemma}   and the  proof of Theorem~\ref{pcubeexample},    $kP$  has a PBW basis   consisting of monomials in the $\tilde x$, $\tilde y $ and $\tilde z $.  Moreover, a  PBW basis element $\tilde x^i\tilde y^j\tilde z^k$ 
of $kP$ with $0\le i,j,k < p$, is an eigenvector in the
$(i+k,j+k)$ eigenspace, where $i+k$ and $j+k$ are
read modulo $p-1$.

\begin{lemma}\label{le:Zp}
We have $\tilde z^p=0$.
\end{lemma}
\begin{proof}
The element $\tilde z^p$ is an eigenvector in the $(1,1)$ eigenspace.  Further,
   $\tilde z^p $   has image $ z^p  =0  \in     \gr_{2p} (kP) $,    hence $\tilde z^p $ is in $J^{2p+1}(kP)$. The PBW basis elements 
in this eigenspace have $i+k$ and $j+k$ congruent to one
modulo $p-1$ and at most $2p-2$, and hence at most $p$,
but then $i+j+2k\le 2p$, so the basis element is not in $J^{2p+1}(kP)$.
It follows that the $(1,1)$ eigenspace in $J^{2p}(kP+1)$ is zero and
so $\tilde z^p=0$.
\end{proof}

\begin{lemma}\label{le:[X,Z]}
The element $[\tilde x,\tilde z]$ is a linear combination of the elements
$\tilde z\tilde y^{p-2}\tilde z \in J^{p+2}(kP)$ and 
$\tilde x^i\tilde y^{2i-1}\tilde x^i\tilde z^{p+1-2i}
\in J^{2p+1}(kP)$ with $1\le i \le (p-1)/2$. Similarly,
$[\tilde y,\tilde z]$ is a linear combination of the elements
$\tilde z\tilde x^{p-2}\tilde z\in J^{p+2}(kP)$ and 
$\tilde y^i\tilde x^{2i-1}\tilde y^i\tilde z^{p+1-2i}
\in J^{2p+1}(kP)$ with $1\le i\le (p-1)/2$.
\end{lemma}
\begin{proof}
We prove the first statement. The proof of the second
is identical, with the roles of $\tilde x$ and $\tilde y$ reversed.

The element $[\tilde x,\tilde z]$ has image $[x,z]=0$ in $\gr_3(kP)$,
and hence lies in $J^4(kP)$. It is in the $(2,1)$ eigenspace,
so we start by identifying the PBW basis elements of $J^4(kP)$
in this eigenspace. These are $\tilde y^{p-2}\tilde z^2$ and 
$\tilde x\tilde y^{p-1}\tilde z\in
J^{p+2}(kP)$ and $\tilde x^{i+1}\tilde y^i\tilde z^{p-i} \in J^{2p+1}(kP)$ with
$1\le i \le p-2$.

However, we also need to make use of symmetry.  Let $\sigma $  be the composition of   the automorphism   of $kP$  which  inverts $g$ and $h$  (and  hence   fixes  $c$)  with  the anti-automorphism  of $kP$  which inverts all elements of $P$. Then $\sigma$  fixes $\tilde x$ and
$\tilde y$, reverses multiplication in $kP$, and negates $\tilde z$. The point
is that $[\tilde x,\tilde z]=\tilde x^2\tilde y
-2\tilde x\tilde y\tilde x+\tilde y\tilde x^2$ is fixed by $\sigma$, whereas
$\sigma$ does not fix all elements of the $(2,1)$ eigenspace.
With this in mind, we modify the PBW basis of this eigenspace
so that the action of $\sigma$ is more transparent. 

The element $\tilde y^{p-2}\tilde z^2$, for example, is not fixed by $\sigma$,
even though it's fixed modulo $J^{p+3}(kP)$. So instead, we
use the element $\tilde z\tilde y^{p-2}\tilde z$, which is equivalent to it modulo
$J^{p+3}(kP)$, and therefore just as good as part of a PBW
basis of $kP$, but is fixed by $\sigma$. Since 
$\sigma(\tilde x\tilde y^{p-1}\tilde z)\equiv 
-\tilde x\tilde y^{p-1}\tilde z-\tilde y^{p-2}\tilde z^2\pmod{J^{p+3}(kP)}$, 
the element $\tilde x\tilde y^{p-1}\tilde z$ is not involved in the
expression for $[\tilde x,\tilde z]$.
So $[\tilde x,\tilde z]$ is congruent to a multiple of 
$\tilde z\tilde y^{p-2}\tilde z$ modulo
$J^{2p+1}(kP)$.

For the linear span of the elements
$\tilde x^{i+1}\tilde y^i\tilde z^{p-i}$, 
since there are no $(2,1)$ eigenvectors lower
in the radical series, reordering the terms in a monomial
has the same effect as in $\cU\Jen_*(kP)$. So we can choose
a basis consisting of the elements 
$\tilde x^i\tilde y^{2i-1}\tilde x^i\tilde z^{p+1-2i}$
$(1\le i \le (p-1)/2$) and the elements 
$\tilde y^i\tilde x^{2i+1}\tilde y^i\tilde z^{p-2i}$
$(1\le i \le (p-3)/2$). The former are $+1$ eigenvectors of
$\sigma$, while the latter are $-1$ eigenvectors. So the
expression for $[\tilde x,\tilde z]$ only involves the former.
\end{proof}

By Lemma~\ref{le:[X,Z]}, we can write
\begin{align}
[\tilde x,\tilde z]&=a_0 \tilde z\tilde y^{p-2}\tilde z 
+ a_1\tilde x\tilde y\tilde x\tilde z^{p-1}
+a_2\tilde x^2\tilde y^3\tilde x^2\tilde z^{p-3}
+\cdots
+a_{(p-1)/2}\tilde x^{\frac{p-1}{2}}\tilde y^{p-2}\tilde x^{\frac{p-1}{2}}\tilde z^2,
\label{eq:XZ} \\
[\tilde y,\tilde z]&=-a_0 \tilde z\tilde x^{p-2}\tilde z 
-a_1\tilde y\tilde x\tilde y\tilde z^{p-1} 
-a_2\tilde y^2\tilde x^3\tilde y^2\tilde z^{p-3}
-\cdots-
a_{(p-1)/2}\tilde x^{\frac{p-1}{2}}\tilde y^{p-2}\tilde x^{\frac{p-1}{2}}\tilde z^2.\label{eq:YZ}
\end{align}
Here, we have used the symmetry of $kP$ which swaps $\tilde x$ and 
$\tilde y$, and negates $\tilde z$, to compare the coefficients in~\eqref{eq:XZ}
and those in~\eqref{eq:YZ}.

\begin{rk}
With the aid of the computer algebra system 
{\sc Magma}~\cite{Bosma/Cannon/Playoust:1997a}  we have determined   the relation~\eqref{eq:XZ}    for small $p$  as follows:
\begin{align*}
p=3:\qquad [\tilde x,\tilde z]&=\tilde z\tilde y\tilde z, \\
p=5:\qquad [\tilde x,\tilde z]&=\tilde z\tilde y^3\tilde z
+2\tilde x\tilde y\tilde x\tilde z^4, \\
p=7:\qquad [\tilde x,\tilde z]&=\tilde z\tilde y^5\tilde z
+4\tilde x\tilde y\tilde x\tilde z^6
+2\tilde x^2\tilde y^3\tilde x^2\tilde z^4. 
\end{align*}
One might surmise that $a_0=1$ and $a_{(p-1)/2}=0$, but we have 
not proved that. Nor have we spotted the general pattern of 
the coefficients. 
\end{rk}

\begin{theorem}
A presentation for $kP$ is given by generators $\tilde x$, $\tilde y$,
$\tilde z$ with the relations~\eqref{eq:XZ}
and~\eqref{eq:YZ} together with
\[ \tilde x^p=\tilde y^p=\tilde z^p=0,\qquad 
[\tilde x,\tilde y]=\tilde z, \]
and relations saying that all words of length at least $4p-3$ in
$\tilde x$ and $\tilde y$ are
equal to zero.
\end{theorem}

\begin{proof}
These relations hold in $kP$ by Lemmas~\ref{le:Zp} and~\ref{le:[X,Z]},
and the fact that $J^{4p-3}(kP)=0$.
Let $\bfA$ be the algebra defined by these generators 
and relations. Then we have a surjective map $\bfA\to kP$
taking $\tilde x$, $\tilde y$ and $\tilde z$ to the elements with the same names.
This induces a map $\gr_*\bfA \to \gr_* kP$. The
relations~\eqref{eq:XZ} and~\eqref{eq:YZ} imply that 
the images $x$, $y$ and $z$ in $\gr_*\bfA$ of $\tilde x$, $\tilde y$ and
$\tilde z$ in  $\bfA$
satisfy $[x,z]=0$ and $[y,z]=0$. Thus
all the relations in $\cU\Jen_*(P)$ hold in $\gr_*\bfA$,
and $\gr_*\bfA\to \gr_* kP$ is an isomorphism. Since
the radical of $\bfA$ is nilpotent, this implies that
$\bfA\to kP$ is an isomorphism.
\end{proof}

Recall from  the proof of Theorem~\ref{pcubeexample}, that  
 setting  $\x=txe$, $\y=sye$ and $\z=stze$   in $\gr_*(kGe)$,
we have that the algebra $\fA$   is  generated by $\x$, $\y$ and $\z$ centralises
$\fM$ in $\cU\Jen_*(P)\rtimes kH$.
Further, these elements satisfy the relations
\[ \x^p=0,\quad \y^p=0,\quad \x\y+\y\x=\z,\quad
\x\z+\z\x=0,\quad \y\z+\z\y=0 \]
(and these imply that $\z^p=0$).

In $kGe$, we set $\tilde \x=t\tilde xe$, $\tilde \y=s\tilde ye$
and $\tilde \z=st\tilde ze$.
The algebra $\tilde\fA$ generated by 
$\tilde\x$, $\tilde\y$ and $\tilde\z$  centralises 
$\fM$ in $kGe$. These elements satisfy
the relations
\[ \tilde \x^p=\tilde \y^p=\tilde \z^p=0,\qquad 
\tilde \x\tilde \y+\tilde \y \tilde \x = \tilde \z \]
together with the following quantised versions of~\eqref{eq:XZ}
and~\eqref{eq:YZ}
\begin{align*}
\tilde \x\tilde \z + \tilde \z \tilde \x & 
=(-1)^{\frac{p-1}{2}}\!
\bigl(a_0 \tilde \z\tilde \y^{p-2}\tilde \z -
a_1\tilde\x\tilde\y\tilde\x\tilde\z^{p-1}-
a_2\tilde\x^2\tilde\y^3\tilde\x^2\tilde\z^{p-3}
-\cdots-
a_{(p-1)/2}\tilde\x^{\frac{p-1}{2}}\tilde\y^{p-2}\tilde\x^{\frac{p-1}{2}}\tilde\z^2\bigr), \\
\tilde\y \tilde\z + \tilde\z \tilde\y &
=(-1)^{\frac{p-1}{2}}\!
\bigl(a_0\tilde\z\tilde\x^{p-2}\tilde\z 
- a_1\tilde\y\tilde\x\tilde\y\tilde\z^{p-1} 
- a_2\tilde\y^2\tilde\x^3\tilde\y^2\tilde\z^{p-3}
-\cdots-
a_{(p-1)/2}\tilde\y^{\frac{p-1}{2}}\tilde\x^{p-2}\tilde\y^{\frac{p-1}{2}}\tilde\z^2\bigr),
\end{align*}
together with relations saying that all 
words of length at least $4p-3$ in $\tilde \x$ and
$\tilde \y$ are equal to zero.\bigskip

Using  {\sc Magma}~\cite{Bosma/Cannon/Playoust:1997a},
in the case $p=3$ we have 
succeeded in finding a short presentation for $kP$
in terms of the generators $\tilde x$ and $\tilde y$. 
In this case, we have $\tilde x=g^{-1}-g$ and $\tilde y=h^{-1}-h$. 
Defining $\tilde z=[\tilde x,\tilde y]$, 
the following relations hold in $kP$.
\begin{align}
\tilde x^3=&0,\qquad \tilde y^3=0,\qquad 
[\tilde x,\tilde y]=\tilde z,\qquad 
[\tilde x,\tilde z]=\tilde z\tilde y\tilde z,\qquad
[\tilde y,\tilde z]=-\tilde z\tilde x\tilde z. \label{eq:kP}
\end{align}
It follows from these relations
that $\tilde z^3=0$, so it is not necessary to include this in the
relations,
and hence that the algebra defined by these relations has
dimension $27$, and is isomorphic to $kP$. 
This is the content of the next theorem. Note, however, that
the proof is difficult, so for some purposes it is
better to adjoin $\tilde z^3=0$ to the above presentation.
We restate and prove the first part of Theorem \ref{27example}.

\begin{theorem}\label{th:27}
Suppose that $p=3$.
The generators $\tilde x$, $\tilde y$ and $\tilde z$ and the 
relations~\eqref{eq:kP} give a presentation for $kP$.
\end{theorem}

\begin{proof}
Since the given elements of $kP$ satisfy these relations, 
it suffices to prove that the algebra defined by the relations
has dimension at most $27$. The crucial point is to prove 
that $\tilde z^3=0$.

It is more convenient to extend the field so that it has a square root
of $-1$, which we denote $\ii$. Then we set 
$a=\tilde x+\ii\tilde y$, $b=\tilde x-\ii\tilde y$,
$c=\ii\tilde z$, and the presentation becomes 
\[ a^3=[b,c]=cbc,\qquad b^3=-[a,c]=cac,
\qquad [a,b]=c, \]
and we must show that $c^3=0$.

We have 
\begin{align*}
(1+c)a(1-c)&=a,\\
(1-c)b(1+c)&=b,
\end{align*}
and so 
\[ (1+c)ab=(1+c)a(1-c)b(1+c)=ab(1+c). \]
Therefore $c$ commutes with $ab$ and with $ba$.

Next,
\[ cab=acb+cacb=abc-acbc+cacb, \] 
and since $cab=abc$ it follows that $c$ commutes with $acb$.
Thus we have
\begin{equation}\label{eq:a4b4} 
cbca=a^4=acbc=cacb=b^4=bcac. 
\end{equation}
 
Since we are in characteristic three, we also have 
\[ [a,[a,c]]=-[a,cac]=-[a,c]ac-c[a,a]c-ca[a,c]=cacac+cacac=-cacac \]
and so 
\[ [a^3,c]=[a,[a,[a,c]]]=-[a,cacac]=-[a,c]acac-ca[a,c]ac-caca[a,c]=-3cacacac=0. \]
Thus $c$ also commutes with $a^3$:
\begin{equation}\label{eq:ca^3}
a^3c=ca^3.
\end{equation}

Next, using~\eqref{eq:a4b4} we have
\begin{align} 
c^3&=cabc-cbac=acbc+cacbc-cbca+cbcac=a^4+b^4c-b^4+a^4c\label{eq:c3}\\
&=-a^4c = -b^4c = -cacbc = -cbcac= -ca^4 = -cb^4.\notag 
\end{align}
Using~\eqref{eq:ca^3} and~\eqref{eq:c3}, we have
\[ a^4c=ca^4 =a^3ca=a^4c+a^3cac=a^4c+ca^4c=a^4c-c^4 \]
and so $c^4=0$.

The fact that $c^4=0$ enables us to write
\begin{align} 
ca&=ac+cac=ac+ac^2+cac^2=ac+ac^2+ac^3+cac^3=a(c+c^2+c^3)\label{eq:ca} \\
cb&=bc-cbc=bc-bc^2+cbc^2=bc-bc^2+bc^3-cbc^3=b(c-c^2+c^3). \label{eq:cb}
\end{align}
We can use these to move copies of $c$ to the end of expressions, at
the expense of accumulating higher powers of $c$. 
Applying this to $a^6=cbc^2bc$, we get $b^2(c^4+{}$higher powers of
$c)$, so we get $a^6=0$. Similarly, we get $b^6=0$. Then when we
do the same with $(ab)^9$, moving $b$ past $a$ using $ba=ab-c$, we get 
\[  a^9b^9+a^8b^8cf_1(c)+a^7b^7c^2f_2(c)+a^6b^6c^3f_3(c)+a^5b^5c^4f_4(c)
  + \dots \]
for suitable polynomials $f_i(c)$. Since $a^6=b^6=c^4=0$,
every term here is zero, and so $(ab)^9=0$. 

Now using~\eqref{eq:c3} and the same method, we have
\[ c^3=-cacbc=-cab(c-c^2+c^3)c=-abc(c-c^2+c^3)c=-abc^3, \]
and so $(1+ab)c^3=0$. Since $ab$ is nilpotent, $(1+ab)$ is invertible,
so this implies that $c^3=0$.

The original relations together with~\eqref{eq:ca},
\eqref{eq:cb} and $c^3=0$ allow us to rewrite every element as a linear
combination of the elements $a^ib^jc^k$ with $0\le i,j,k < 3$,
so the algebra has dimension at most $27$, and we are done.
\end{proof}

\begin{proof}[Proof of Theorem \ref{27example}]
The relations for $kP$ are proved in Theorem \ref{th:27}.
As above, using $\tilde x=g^{-1}-g$, $\tilde y = h^{-1}-h$, 
and $\tilde z=[\tilde x,\tilde y]$, 
to obtain generators for the basic algebra for $kGe$, we set 
$\tilde \x = t\tilde xe$ and $\tilde \y = s\tilde ye$, 
$\tilde \z = st\tilde ze$.
These satisfy
\[
\tilde \x^3=0, \qquad \tilde \y^3=0,\qquad 
\tilde \x\tilde \y+\tilde \y\tilde \x=\tilde \z, \qquad
\tilde \x \tilde \z +
\tilde \z \tilde \x =
 -\tilde \z\tilde \y\tilde \z, \qquad
\tilde \y \tilde \z + \tilde \z \tilde \y=
-\tilde \z \tilde \x \tilde \z.
\]
Furthermore, the algebra defined by these relations again has
dimension $27$, and is hence isomorphic to the basic algebra of $kGe$.
\end{proof}

\begin{rk}
The above relations for $kGe$  are 
a quantised version of the relations for $kP$. These relations imply
that $\tilde \z^3=0$, and adjoining this relation makes the
presentation easier to work with if desired.

These presentations can be lifted to give integral presentations. The algebra
$\cO P$ has a presentation with corresponding generators $\hat x$, 
$\hat y$ and $\hat z$ 
subject to
\[ \hat x^3+3\hat x=0,\qquad \hat y^3+3\hat y=0,
\qquad [\hat x,\hat y]=\hat z,\qquad
2[\hat x,\hat z]=-\hat z\hat y\hat z,
\qquad 2[\hat y,\hat z]=\hat z\hat x\hat z, \]
while $\cO Ge$ is generated by $\hat \x$, $\hat \y$ and $\hat \z$ 
subject to
\[ \hat \x^3=3\hat \x,\qquad \hat \y^3 =3\hat \y,\qquad
\hat \x\hat \y+\hat \y\hat \x = \hat \z,\qquad
2\hat \x \hat \z + 2\hat \z \hat \x = \hat \z\hat \y\hat \z,\qquad 
2\hat \y \hat \z + 2\hat \z \hat \y =  \hat \z \hat \x \hat \z. \]
The expressions for $\hat z^3$ in $\cO P$ and for 
$\hat\z^3$ in $\cO Ge$, lifting the fact that they cube to 
zero modulo three, are ugly even though they follow from the
presentations above.
\end{rk}

\section{\texorpdfstring{Example: $2^{1+4}{:\,}3^{1+2}$ in
    characteristic two}{Example: 2¹ᐩ⁴:3¹ᐩ² in characteristic two}}
\label{se:char2}

The examples in the last section were at odd primes
for extraspecial groups of order $p^3$. In this section
we give an example in characteristic two with an extraspecial group of
order $2^5$.

Let $P$ be an extraspecial group $2^{1+4}$ which is a 
central product of two copies of the quaternion
group of order eight, and let $H$ be an extraspecial
group $3^{1+2}$ of exponent three. We let the centre $Z\cong\bZ/3$
of $H$ act trivially on $P$, and the elementary abelian quotient
act as the automorphisms of order three on the two quaternion
central factors of $P$, and we set $G=P\rtimes H$. Thus the quotient
$G/Z\cong SL(2,3)\circ SL(2,3)$ is a central product of two copies of the 
group $SL(2,3)$ of order $24$.

More precisely, we let 
\begin{align*} 
P &= \langle g_1,g_2, h_1,h_2, c \mid 
g_1^2=h_1^2=[g_1,h_1]=g_2^2=h_2^2=[g_2,h_2]=c,\\ 
&\qquad [g_1,c]=[g_2,c]=[h_1,c]=[h_2,c]=
[g_1,g_2]=[g_1,h_2]=[h_1,g_2]=[h_1,h_2]=c^2=1 \rangle, \\
H &= \langle s_1,s_2, t \mid
s_1^3=s_2^3=1,\ [s_1,s_2]=t,\ [s_1,t]=[s_2,t]=t^3=1\rangle. 
\end{align*}
Let $H$ act on $P$ with $Z=\langle t\rangle$ acting trivially, and
\begin{gather*}
s_1g_1s_1^{-1}=h_1, \qquad 
s_1h_1s_1^{-1}=g_1h_1,\qquad 
s_1g_2s_1^{-1}=g_2, \qquad 
s_1h_2s_1^{-1}=h_2,\qquad
s_1cs_1^{-1}=c, \\
s_2g_1s_2^{-1}=g_1, \qquad 
s_2h_1s_2^{-1}=h_1,\qquad
s_2g_2s_2^{-1}=h_2, \qquad
s_2h_2s_2^{-1}=g_2h_2, \qquad 
s_2cs_2^{-1}=c.
\end{gather*}
Let $k$ be a field of characteristic two containing 
$\bF_4=\{0,1,\omega,\bar\omega\}$. A basis of eigenvectors
in $\gr_1(kP)$ is given by
\[ x_i = \bar\omega(g_i-1)+\omega(h_i-1),\qquad
y_i=\omega(g_i-1)+\bar\omega(h_i-1)
\qquad (i=1,2). \]
These give the following presentation for $\gr_*(kP)$.
\[ x_i^2=0,\quad y_i^2=0,\quad
  [x_1,x_2]=[y_1,y_2]=[x_1,y_2]=[x_2,y_1]=0,
\quad [x_1,y_1]=[x_2,y_2]. \]
A lift of $x_i$ and $y_i$ to eigenvectors 
complementing $J^2(kP)$ in $J(kP)$ is given by the elements
\[ \tilde x_i = \omega g_i + \bar\omega h_i +  g_ih_i, \qquad  
\tilde y_i = \bar\omega g_i + \omega h_i +  g_ih_i\qquad
(i =1,2). \]
The relations lift to
\begin{gather*} 
\tilde x_i^2=\tilde y_i\tilde x_i\tilde y_i,
\qquad \tilde y_i^2=\tilde x_i\tilde y_i\tilde x_i,\qquad 
\tilde x_i^4=0 \qquad (i=1,2),\\
[\tilde x_1,\tilde x_2]=[\tilde y_1,\tilde y_2]=
[\tilde x_1,\tilde y_2]=[\tilde x_2,\tilde y_1]=0, \qquad
[\tilde x_1,\tilde y_1]+\tilde x_1^3=[\tilde x_2,\tilde y_2]+\tilde x_2^3
\end{gather*}
(both sides in the last relation are equal to 
$(1+c)$). 
Using the  radical filtration,  it is not hard to check that these relations define a $k$-algebra of 
dimension  at most  $32$, which is therefore isomorphic to $kP$.

The action of $H$ on $kP$ with respect to these generators is given by
\[ g\tilde x_ig^{-1}=\psi_i(g)^{-1}\tilde x_i,\qquad 
g\tilde y_ig^{-1}=\psi_i(g)\tilde y_i\qquad  (g\in H) \]
where $\psi_1(s_1)=\psi_2(s_2)=\omega$, $\psi_1(s_2)=\psi_2(s_1)=1$.

Set 
\[ e_0=1+t+t^2,\qquad e=1+\bar\omega t + \omega t^2,
\qquad \bar e = 1 + \omega t +\bar\omega t^2. \] 
Then $kG$ has three blocks, the principal block $kGe_0$, and
two non-principal blocks $kGe$ and $kG\bar e$.
We examine the non-principal block $kGe$, the other is similar.

We have 
$ s_1s_2e=\omega s_2s_1e$.  
We set
\[
\tilde \x_1 = s_2\tilde x_1e,\qquad 
\tilde \x_2 = s_1^{-1}\tilde x_2e, \qquad
\tilde \y_1 = s_2^{-1}\tilde y_1e,\qquad
\tilde \y_2 = s_1\tilde y_2e.
\]
These commute with $\fM=kHe$, and generate the
subalgebra $\fA$, so that $kGe\cong \Mat_{3}(\fA)$. 
They satisfy the relations:
\begin{gather*} 
\tilde \x_i^2=\tilde \y_i\tilde \x_i\tilde \y_i,\qquad
\tilde \y_i^2=\tilde \x_i\tilde \y_i\tilde \x_i,\qquad
\tilde \x_i^4=0 \qquad (i=1,2), \\
\tilde \x_1\tilde \x_2=
\bar\omega\tilde \x_2\tilde \x_1, \qquad
\tilde \x_1\tilde \y_2=
\omega \tilde \y_2\tilde \x_1, \\
\tilde \y_1\tilde \x_2=
\omega \tilde \x_2\tilde \y_1, \qquad\ 
\tilde \y_1\tilde \y_2=
\bar \omega\tilde \y_2\tilde \y_1, \\
[\tilde \x_1,\tilde \y_1]+\tilde \x_1^3=
[\tilde \x_2,\tilde \y_2]+\tilde \x_2^3
\end{gather*}
(both sides in the last relation are equal to 
$(1+c)e$). These are identical to the relations for $kP$ apart from
the commutation relations, which have been quantised by the
introduction of factors $\omega$ and $\bar\omega$.

\section{Appendix: Errata}\label{se:errata} 

The present paper supersedes most of our 
previous paper~\cite{Benson/Kessar/Linckelmann:2019a}. In that paper
there are a number of minor errors, mostly in the calculations in Section~4,
which have been corrected in the present work. 
We give a list of those errors in \cite{Benson/Kessar/Linckelmann:2019a}. 

In the statements of Theorem 1.2 and Corollary 1.3, it should read `...quantised version
of $k(P\rtimes Z(H)/Z)$' (and not `...of $k(P\rtimes L)$.')

In the 3rd line of the proof of Proposition 3.1 insert the word `abelian' between
`maximal' and `subgroup' (as is done correctly in line 2 and line 4 of that proof).

On page 1441, third line from the bottom, insert \emph{faithful}:
\begin{center}
``\dots and a faithful linear
character $\chi\colon Z \to k^\times$\dots''\bigskip
\end{center}

On page 1443, in the first line, $\rho(g)\colon h \mapsto
\chi([h,g])$. The display on the third line should read
\[ \bar\rho\colon H/Z(H) \to \mathsf{Hom}(H/Z(H),k^\times) \]

On page 1444, line four should begin ``where $\rho(g_{i,\phi})(h)=\chi([h,g_{i,\phi}])$.''
The third line of the proof of Lemma~4.8 should begin with
$eg_{i,\phi}h=e\chi([h,g_{i,\phi}])^{-1}hg_{i,\phi}$.
The displayed equation on the fourth line of the
proof of Lemma~4.8 should read
\[ \xi_\phi(h)^{-1}(g_{i,\phi}w_i)(e_\phi\cdot h) = 
\xi_\phi(h)^{-1}\psi_i(h)^{-1}\chi([h,g_{i,\phi}])^{-1}(e_{\phi\psi_i}h)(g_{i,\phi}w_i). \] 
In Definition 4.9 and the four lines following, $kHe$
should be $k\tilde Ge$ four times. The dimension of $\mathfrak A$
should be given as $|P|\cdot |Z(H):Z|$ and not $|P|\cdot |H:Z(H)|$.\bigskip

On page 1445, in Lemma~4.12\,(2), in the displayed equation the last
$w_i$ should be $w_j$. The scalar $q_{i,j,\phi}$ should equal
$\psi_i(g_{j,\phi}z_{i,j,\phi})\psi_j(g^{-1}_{i,\phi}z_{i,j,\phi})\phi(z_{i,j,\phi})$ rather
than
$\phi(z_{i,j,\phi})$. Similarly, in
Lemma~4.12\,(3) the scalar $q_{i,j,\phi}$ should equal 
$\psi_i(g_{j,\phi})\psi_j(g^{-1}_{i,\phi})\chi(z_{i,j,\phi})$ rather than $\chi(z_{i,j,\phi})$.
In the second line of the proof of Lemma~4.12\,(1), $g_{j,\phi\phi_i}$
should be $g_{j,\phi\psi_i}$.
The computation that was suppressed in the proof of Lemma~4.12\,(2) 
uses~(1) and equation~(4.7). It
is similar to the computation in Theorem~\ref{th:rels}
above, which we have spelled out in detail.
There is a missing $Z$ in the third to last line of the proof of Lemma~4.12\,(3), and
$g_ji$ should be $g_i$ in the second to last line.\bigskip

On page 1447, in Theorem 4.15 and Corollary 4.16, $kGe$ should
be $k\tilde Ge$ five times.

\bibliographystyle{amsplain}
\bibliography{../repcoh}

\end{document}